\documentclass{amsart}

\usepackage{xypic}
 \input xy
\xyoption{all}
\usepackage{epsfig}
\usepackage{amsthm}
\usepackage{amssymb}
\usepackage{amsmath}
\usepackage{amscd}

\newcommand{\bg}{\begin{equation}}
\newcommand{\ed}{\end{equation}}
\newcommand{\bga}{\begin{eqnarray}}
\newcommand{\eda}{\end{eqnarray}}
\newcommand{\pf}{\textbf{Proof:\ }}

\def\cbdu{\par{\raggedleft$\Box$\par}}

\newtheorem {Theorem}  {Theorem}

\numberwithin{Theorem}{section}

\newtheorem {Lemma}[Theorem]  {Lemma}
\newtheorem {Proposition}[Theorem]{Proposition}
\theoremstyle{definition}

\theoremstyle{remark}
\newtheorem{Remark}[Theorem]{\bf Remark}

\newtheorem {Corollary}[Theorem]{\bf Corollary}

\expandafter\chardef\csname pre amssym.def
at\endcsname=\the\catcode`\@ \catcode`\@=11
\def\undefine#1{\let#1\undefined}
\def\newsymbol#1#2#3#4#5{\let\next@\relax
 \ifnum#2=\@ne\let\next@\msafam@\else
 \ifnum#2=\tw@\let\next@\msbfam@\fi\fi
 \mathchardef#1="#3\next@#4#5}
\def\mathhexbox@#1#2#3{\relax
 \ifmmode\mathpalette{}{\m@th\mathchar"#1#2#3}%
 \else\leavevmode\hbox{$\m@th\mathchar"#1#2#3$}\fi}
\def\hexnumber@#1{\ifcase#1 0\or 1\or 2\or 3\or 4\or 5\or 6\or 7\or 8\or
 9\or A\or B\or C\or D\or E\or F\fi}

\font\teneufm=eufm10 \font\seveneufm=eufm7 \font\fiveeufm=eufm5
\newfam\eufmfam
\textfont\eufmfam=\teneufm \scriptfont\eufmfam=\seveneufm
\scriptscriptfont\eufmfam=\fiveeufm

\catcode`\@=\csname pre amssym.def at\endcsname

\newcounter{remark}
\usepackage{color}

\newcommand{\e}{\epsilon}

\newcommand{\R}{\mathbf{R}}

\def  \R   {{\mathbb R}}

\def  \12  {{\frac{1}{2}}}

\def\build#1_#2^#3{\mathrel{\mathop{\kern 0pt#1}\limits_{#2}^{#3}}}

\begin{document}

\title[Asymptotic Behavior of LCD in $H^m(\R^3)$]{Asymptotic Behavior of Solutions to the Liquid Crystal Systems in $H^m(\mathbb{R}^3)$}

\author[Mimi Dai]{ Mimi Dai}
\address{Department of Mathematics, UC Santa Cruz, Santa Cruz, CA 95064,USA}
\email{mdai@ucsc.edu}
\author[Maria Schonbek] {Maria Schonbek}
\address{Department of MAthematics, UC Santa Cruz, Santa Cruz, CA 95064, USA}
\email{schonbek@ucsc.edu}


\thanks{The work of Mimi Dai was partially supported by NSF Grant DMS-0900909;}
\thanks{The work of
M. Schonbek was partially supported by NSF Grant DMS-0900909}


\begin{abstract}
 In this paper we study the large time behavior of  regular solutions to a nematic liquid crystals system in Sobolev spaces $H^m(\R^3)$ for $m\geq 0$.We obtain optimal decay rates in $H^m(R^3)$  spaces, in the sense that the rates coincide with the rates of the underlying linear counterpart.  The fluid under consideration has  constant  density  and small  initial data.  \end{abstract}

\maketitle

\section{Introduction}
In this paper we consider the asymptotic behavior of solutions to the simplified model of nematic liquid crystals (LCD) with constant density in Sobolev spaces $H^m(\R^3)$ for $m\geq 0$:

\begin{equation}\begin{split}\label{LCD}
u_t+u\cdot\nabla u+\nabla \pi =\nu\triangle u -\nabla\cdot(\nabla d\otimes\nabla d),\\
d_t+u\cdot\nabla d =\triangle d-f(d),\\
\nabla\cdot u =0.
\end{split}
\end{equation}
The equations are considered in $\R^3\times(0, T)$. Here,
 $\pi: \R^3\times[0,T]\to\mathbb{R}$ is the fluid pressure,
 $u: \R^3\times[0,T]\to\mathbb{R}^3$  the fluid velocity,
and $d: \R^3\times[0,T]\to\mathbb{R}^3$  the director field representing the alignment of the molecules.
The constant $\nu >0$ stands for the viscosity coefficient.  Without loss of generality, by scaling, we set $\nu =1$.  The forcing term $\nabla d\otimes\nabla d$ in the equation of  conservation of momentum denotes the $3\times 3$ matrix whose $ij$-th entry is given by $``\nabla_i d\cdot\nabla_j d", 1\leq i,j\leq 3$. This force $\nabla d\otimes\nabla d$ is the stress tensor of the energy about the director field $d$, where the   energy is given by:
$$
\frac 12 \int_{\R^3} |\nabla d|^2 dx + \int_{\R^3}F(d)dx
$$
where
$$
F(d)=\frac{1}{4\eta^2}(|d|^2-1)^2, \quad f(d) = \nabla F(d) = \frac{1}{\eta^2}(|d|^2-1)d,
$$
for a constant $\eta$ in this paper. The  $F(d)$ is the penalty term of the Ginzburg-Landau approximation of the original free energy for the director field with unit length.\\
We consider the following  initial conditions:
\bg\label{initu}
u(x,0) =u_0(x), \ \ \ \nabla\cdot u_0 =0,
\ed
\bg\label{initd}
d(x,0) =d_0(x), \ \ \ |d_0(x)|= 1,
\ed
and
\bg\label{bd}
u_0\in H^m(\R^3), \ \ d_0-w_0\in H^{m+1}(\R^3),
\ed
for any integer $m\geq 1$ with a fixed vector $w_0\in S^2$, that is, $|w_0|=1$.

The flow of nematic liquid crystals can be treated as slow moving particles where the fluid velocity and the alignment of the particles influence each other. The hydrodynamic theory of liquid crystals  was established by Ericksen \cite{Er0, Er1} and Leslie \cite{Le0, Le1} in the 1960's.
As F.M. Leslie points out in his 1968 paper: ``liquid crystals are states of matter which are capable of flow, and in which the molecular arrangements  give rise to a preferred direction".  There is a vast literature on the hydrodynamic of  liquid crystal systems. For background we list  a few, with no intention to be complete:  \cite{EK, HKL, Kin, LL, LL2, LL1, Cal, CC, CC1, CDLL, Wu, Liu, JT}.  In particular, the asymptotic behavior of solutions to the flow of nematic liquid crystals was studied for bounded domains in \cite{LL, Wu}.  It was shown in \cite{Wu} that, with suitable initial conditions, the velocity converges to zero and the direction field converges to the steady solution
to the following equation
\begin{equation}\label{eq:steady}
\begin{cases}
-\Delta d+f(d)=0, x\in\Omega\\
d(x)=d_0(x), x\in\partial\Omega.
\end{cases}
\end{equation}
In  \cite{Wu}, Lemma 2.1 the \L ojasiewicz-Simon inequality  is used to derive the convergence when $\Omega$ is a bounded domain.

\noindent In our previous work \cite{DQS}, we established a preliminary  decay rate for the solutions in $\R^3$  to \eqref{LCD},
subject to the additional condition on  the director field which insures that the  initial director  field tends to a constant unit vector $w_0$,  as the space variable tends to infinity:
\begin{equation}\label{d0w0}
\lim_{|x|\to\infty} d_0(x) = w_0.
\end{equation}
This behavior at infinity of the initial  director field allows  to obtain the stability without needing
the Liapunov reduction and \L ojasiewicz-Simon inequality, since $w_0$ is a non-degenerate steady solution to \eqref{eq:steady}.

\medskip

The paper is organized as follows, in section 2 we recall some previous results,  give some preliminary estimates and state the main Theorem. In section 3 we obtain the decay rate for the velocity in $L^2(\R^3)$. Section 4 deals with the decay rate for gradient of the director field in $L^2(\R^3)$. The last section gives  decay rates in $H^m$ for the velocity and the director vector. All the obtained rates are optimal. The rates obtained  improve the rates obtained in \cite{DQS}. We work with  regular solutions with small data.
The main tools used are the Fourier splitting method and appropriate energy estimates.\\

\begin{Remark}We note that for the $L^2$  the decay rates can be obtained  also for weak solutions. In this case one shows the decay for approximations in the form we do in the paper and passing to the limit it will follow for the solutions. We expect that for higher derivatives one can obtain the decay for  approximations, and show that eventually the solutions become small and the use the results in this paper.
\end{Remark}
   
\section{Preliminaries} 
We start from the basic energy estimates and Ladyzhenskaya higher order energy estimates \cite{LS,DQS}. We establish that the velocity $u$ converges to $0$ with decay rate $(1+t)^{-1/4}$ in $L^2(\R^3)$ and that the convergence of the direction field $d$ to the constant steady solution $w_0$ with the decay rate $(1+t)^{-\frac{3}{2}(1-\frac{1}{p})}$ in $L^p(\R^3)$ for any $p>1$. The main ingredients are the Fourier splitting method and the Gagliardo-Nirenberg interpolation techniques. \\

We recall that for small data there exists a smooth solution given by

\begin{Theorem}\cite{DQS}\label{exist} Let $u_0$ and $d_0$ satisfy (\ref{initu})-(\ref{bd}). Assume that $u_0\in H^1(\R^3)$ and $d_0- w_0 \in H^{2}(\R^3)\cap L^1(\R^3)$ for a unit vector $w_0$. There is a positive small number $\e_0$ such that if
\begin{equation}\label{smalldata}
\|u_0\|_{H^1(\R^3)}^2+\|d_0 - w_0\|_{H^{2}(\R^3)}^2\leq\epsilon_0,
\end{equation}
then the system (\ref{LCD}) has a classical solution $(u,\pi,d)$ in the time period $(0, T)$, for all $T>0$. That is, for some $\alpha\in(0,1)$
\begin{equation}\label{reg}\begin{split}
u\in C^{1+\alpha/2,2+\alpha}((0,T)\times\R^3)\\
\nabla \pi\in C^{\alpha/2,\alpha}((0,T)\times\R^3)\\
d\in C^{1+\alpha/2,2+\alpha}((0,T)\times\R^3).
\end{split}\end{equation}
And the solution $(u,\pi,d)$ satisfies the following basic energy estimate and higher order energy estimate 
\begin{align}\label{ex:basic}
&\int_{\R^3}|u|^2+|\nabla d|^2+2F(d)dx+2\int_0^T\int_{\R^3}|\nabla u|^2+|\Delta d-f(d)|^2dxdt\\
&\leq \|u_0\|_{L^2(\R^3)}^2+\|\nabla d_0\|_{L^2(\R^3)}^2\notag
\end{align}
\begin{align}\label{ex:lady}
&\int_{\R^3}|\nabla u|^2+|\Delta d|^2dx+\int_0^T\int_{\R^3}|\Delta u|^2+|\nabla\Delta d|^2dxdt\\
&\leq C(\|u_0\|_{H^1(\R^3)}^2+\|d_0 - w_0\|_{H^{2}(\R^3)}^2).\notag
\end{align}
with the constant $C$ depending only on initial data and on $\eta$.
\end{Theorem}

In \cite{DQS} we established the following decay result for the regular solutions obtained in Theorem \ref{exist}:

\begin{Theorem}\label{Mthm}
Let $(u,\pi,d)$ be a regular solution to the system (\ref{LCD}) with initial data (\ref{initu}) and (\ref{initd}) and boundary condition (\ref{bd}) for $m=1$. Assume additionally $u_0\in L^1(\R^3)$ and $d_0-w_0\in L^p(\R^3)$, for any $1\leq p<\infty$ and a unit vector $w_0$. There exists a small number $\e_0>0$ such that if
\begin{equation}\label{smalldataH2}
\|u_0\|_{H^1(\R^3)}^2+\|d_0-w_0\|_{H^2(\R^3)}^2\leq\epsilon_0,
\end{equation}
then
\bg\label{ddecay}
\|d(\cdot,t)-w_0\|_{L^p(\R^3)}\leq C_p\|d_0-w_0\|_{L^p(\R^3)}(1+t)^{-\frac{3}{2}(1-\frac{1}{p})},
\ed
\bg\label{dh1decay}
\|\nabla (d(\cdot,t)-w_0)\|_{L^2(\R^3)}^2\leq C(1+t)^{-\frac{3}{4}},
\ed
\bg\label{udecay}
\|u(\cdot,t)\|_{L^2(\R^3)}^2\leq C(1+t)^{-\frac{1}{2}},
\ed
where the various constants $C$ depend only on initial data, $C_p$  depends on the data and p.
\end{Theorem}

We first need to establish  energy estimates for the velocity $u$ in $H^m(\R^3)$ and the director field $d$ in $H^{m+1}(\R^3)$, for all $m\geq 1$, starting with more regular initial data. We will use the following notation
\begin{equation}\label{eq:phi}
\Phi^2_k(t)=\|D^k u\|_{L^2(\R^3)}^2+\|D^{k+1} d\|_{L^2(\R^3)}^2.
\end{equation}
\begin{equation}\label{eq:psi}
\Psi^2_m(t)= \sum_{k=0}^m \Phi^2_k(t)
\end{equation}
\begin{Theorem}\label{thm:Hexist}
Let $(u,\pi,d)$ be a regular solution to the system (\ref{LCD}) obtained in Theorem \ref{exist}. Assume (\ref{bd}) holds. Then the solution $(u,\pi,d)$ satisfies, for all $m\geq 1$ 
\begin{align}\label{ex:Hlady}
&\int_{\R^3}\Psi_m^2 dx+\int_0^T\int_{\R^3}\Psi_{m+1}^2dxdt\\
&\leq C_m(\|u_0\|_{H^m(\R^3)}^2+\|d_0 - w_0\|_{H^{m+1}(\R^3)}^2).\notag
\end{align}
The constant $C_m$ depends only on initial data,  $\eta$ and $C_k$ for $k=0,\ldots,m-1$, and $C_0$ depends only on the data.

\end{Theorem}

Our main result reads as follows

\begin{Theorem}\label{Mthm:hm}
Assume the initial data $(u_0,d_0)$ satisfies (\ref{initu})-(\ref{bd}),
 $u_0\in L^1(\R^3)$ and $d-w_0\in L^p(\R^3)$ for $p\geq 1$. Let $(u,\pi,d)$ be the regular solution obtained in Theorem \ref{exist}. 
Then, for all $m\geq 0$
\bg\label{dec:hm-u}
\|u(\cdot,t)\|_{H^m(\R^3)}^2\leq C_m(1+t)^{-(m+\frac{3}{2})},
\ed
\bg\label{dec:hm1-d}
\|d(\cdot,t)-w_0\|_{H^{m}(\R^3)}^2\leq C_m(1+t)^{-(m+\frac{3}{2})},
\ed
\bg\label{dec:infty-u}
\|D^mu(\cdot,t)\|_{L^\infty(\R^3)}\leq C_m(1+t)^{-\frac{m+3}{2}},
\ed
and
\bg\label{dec:infty-d}
\|D^m(d(\cdot,t)-w_0)\|_{L^\infty(\R^3)}\leq C_m(1+t)^{-\frac{m+3}{2}}.
\ed
Here the various constants $C_m$ depend only on initial data and $m$.
\end{Theorem}

\begin{Remark}
In \cite{DQS} where we established Theorem \ref{Mthm}, we remarked  that   inequality (\ref{ddecay}) doest not include the limit case $p=\infty$ which would  produces the decay of $d-w_0$ in $L^\infty(\R^3)$ with the rate $(1+t)^{-3/2}$. The reason was that the small constant $\epsilon_0$ in Theorem \ref{Mthm} would be forced to be zero if $p=\infty$. In this paper, we are able to obtain the optimal decay rate $(1+t)^{-3/2}$ for $d-w_0$ in $L^\infty(\R^3)$, which is included in (\ref{dec:infty-d}).  
\end{Remark}

  In \cite{Sch94}, the authors studied the large time behavior of solutions to Navier-Stokes equation in $H^m(\R^n)$ for all $n\leq 5$.  The decay estimate (\ref{dec:hm-u}) for the velocity obtained in Theorem \ref{Mthm:hm} coincides with the result in \cite{Sch94} for Navier-Stokes equation.

\bigskip

\section{Higher order energy estimates}

In this section  we  show that solutions with initial data in $H^m(\R^3)$ remain in $H^m(\R^3)$
for all time.  We prove the energy inequality (\ref{ex:Hlady}) in Theorem \ref{thm:Hexist}. The main tool used is a modified Ladyzhenskaya energy method (see \cite{DQSr}).\\

In the sequel  we need  to use a Gagliardo-Nirenberg interpolation inequality. For completeness we recall the inequality here
\begin{Proposition}\label{PPGN}\cite{F}
 Let $w\in W^{m,p}(\R^n)\cap L^q(\R^n )$, for $1\leq p\leq \infty$ and $1\leq q\leq \infty$. Then
\bg\label{GN}
\|D^kw\|_{L^r(\R^n )}\leq C_m\|D^mw\|_{L^p(\R^n)}^a\|w\|_{L^q(\R^n )}^{1-a}
\ed
for any integer $k\in[0,m-1]$, where
\bg\label{parameter}
\frac{1}{r}=\frac{k}{n}+a(\frac{1}{p}-\frac{m}{n})+(1-a)\frac{1}{q}
\ed
with $a\in[\frac{k}{m},1]$, either if $p=1$ or $p>1$ and $m-k-\frac{n}{p}\notin\mathcal{N}\cup\left\{0\right\}$, while $a\in[\frac{k}{m},1)$, if $p>1$ and $m-k-\frac{n}{p}\in\mathcal{N}\cup\left\{0\right\}$.
\end{Proposition}

We now prove Theorem \ref{thm:Hexist}.\\

{\bf{Proof of Theorem \ref{thm:Hexist}:}}
The proof proceeds by induction. We recall that inequality (\ref{ex:Hlady}) for $m=0$ is easily obtained by energy estimates and has been established in \cite{DQS}. Suppose (\ref{ex:Hlady}) holds for $m=1, 2, \ldots, j-1$, that is
\begin{align}\label{ex:hyp}
&\int_{\R^3}\Psi_{j-1}^2 dx+\int_0^T\int_{\R^3}\Psi_{j}^2dxdt\\
&\leq C_{j-1}(\|u_0\|_{H^m(\R^3)}^2+\|d_0 - w_0\|_{H^{m+1}(\R^3)}^2).\notag
\end{align}

 We need show (\ref{ex:Hlady}) holds for $m=j$. Let $0\leq k\leq j$. Recall that
\begin{equation}\label{eq:phi1}
\Phi^2_k(t)=\|D^k u\|_{L^2(\R^3)}^2+\|D^{k+1} d\|_{L^2(\R^3)}^2.
\end{equation}

The idea is to establish an ordinary differential inequality for $\Phi^2_k(t)$ using the equations in system
(\ref{LCD}) and Gagliardo-Nirenberg interpolation inequalities, and then sum all the terms $\Phi^2_0(t), \Phi^2_1(t), \ldots, \Phi^2_k(t)$.  \\
\begin{Remark}
We note that since u, d are in $H^m$ and are regular,there are no boundaries that we have worry about in the subsequent integration by parts.
\end{Remark}
Taking derivative of $\Phi^2_k(t)$ with respect to time yields
\begin{align}\label{eq:dtphi}
\frac{1}{2}\frac{d}{dt}\Phi^2_k(t)&=\int_{\R^3}D^k u\cdot D^k u_tdx+\int_{\R^3}D^{k+1} d\cdot D^{k+1} d_tdx\\
&=-\int_{\R^3}D^{k+1} u\cdot D^{k-1} u_tdx-\int_{\R^3}D^{k+2} d\cdot D^k d_tdx\notag
\end{align}

Taking the $(k-1)$-th spatial derivative on the first equation in (\ref{LCD}) gives
\begin{equation}\notag
D^{k-1} u_t=D^{k-1}\Delta u-D^{k-1}(u\cdot\nabla u)-D^{k-1}\nabla \pi-D^{k-1}(\nabla\cdot(\nabla d\otimes\nabla d)).
\end{equation}
Thus, 
\begin{align}\label{eq:dNSE}
&-\int_{\R^3}D^{k+1}uD^{k-1}u_t\\
&=-\int_{\R^3}|D^{k+1}u|^2dx+\int_{\R^3}D^{k+1}uD^{k-1}(u\cdot\nabla u)dx\notag\\
&+\int_{\R^3}D^{k+1}uD^{k-1}(\nabla\cdot(\nabla d\otimes\nabla d))dx\notag\\
&\leq -\frac{1}{2}\int_{\R^3}|D^{k+1}u|^2dx+4\int_{\R^3}|D^k(u\otimes u)|^2dx+4\int_{\R^3}|D^k(\nabla d\otimes \nabla d)|^2dx\notag.
\end{align}
Taking the $k$-th spatial derivative on the second equation in (\ref{LCD}) gives
\begin{equation}\notag
D^k d_t=D^k\Delta d-D^k(u\cdot\nabla d)-D^k(f(d)).
\end{equation}
Thus,
\begin{align}\label{eq:ddirector}
&-\int_{\R^3}D^{k+2}dD^{k}d_t\\
&=-\int_{\R^3}|D^{k+2}d|^2dx+\int_{\R^3}D^{k+2}dD^{k}(u\cdot\nabla d)dx+\int_{\R^3}D^{k+2}dD^{k}f(d)dx\notag\\
&\leq -\frac{1}{2}\int_{\R^3}|D^{k+2}d|^2dx+4\int_{\R^3}|D^k(u\cdot\nabla d)|^2dx+4\frac{1}{\eta^2}\int_{\R^3}|D^k((|d|^2-1)d)|^2dx\notag.
\end{align}

Inserting (\ref{eq:dNSE}) and (\ref{eq:ddirector}) into (\ref{eq:dtphi}) yields,
\begin{align}\label{eq:dtphi1}
&\frac{d}{dt}\Phi^2(t)+\int_{\R^3}|D^{k+1} u|^2dx+\int_{\R^3}|D^{k+2} d|^2dx\\
&\leq 4\left[\int_{\R^3}|D^k(u\otimes u)|^2dx+\int_{\R^3}|D^k(\nabla d\otimes \nabla d)|^2dx\right]\notag\\
&+4\left[\int_{\R^3}|D^k(u\cdot\nabla d)|^2dx+\frac{1}{\eta^2}\int_{\R^3}|D^k((|d|^2-1)d)|^2dx\right]\notag\\
&\equiv 4[I+II+III+IV]\notag.
\end{align}
The four terms on the right hand side of (\ref{eq:dtphi1}) are estimated as follows.\\
For $I$, we have
\begin{align}\label{ineq:i}
I&\leq 2\|u\|_{L^\infty}^2\|D^ku\|_{L^2}^2+2\|\nabla u\|_{L^\infty}^2\|D^{k-1}u\|_{L^2}^2+I_1,
\end{align}
with
\begin{equation}\label{i1}
I_1=
\begin{cases}
&0, \ \ \mbox { if } k=2,3\\
&\sum_{r=2}^{[(k-2)/2]}\|D^ru\|_{L^\infty}^2\|D^{k-r}u\|_{L^2}^2, \ \ \mbox { if } k\geq 4.
\end{cases}
\end{equation}
Where  $[a]$ denotes the largest integer that is less than or equal $a$.\\
Recall $u$ is regular by (\ref{reg}), $\|u\|_{L^\infty}$ and $\|\nabla u\|_{L^\infty}$ are bounded by some fixed constant $C$. By the induction hypothesis (\ref{ex:hyp}), $\|D^{k-1}u\|_{L^2}^2$ is bounded by $\Psi_k^2$. Thus,
\begin{equation}\label{ineq:imid}
\|u\|_{L^\infty}^2\|D^ku\|_{L^2}^2+\|\nabla u\|_{L^\infty}^2\|D^{k-1}u\|_{L^2}^2\leq C[\Psi^2_k(t)].
\end{equation}
By the remark above  $C$ denotes an absolute constant.
Applying Gagliardo-Nirenberg interpolation inequality,
\begin{equation}\notag
\|D^ru\|_{L^\infty}\leq C\|D^{r+2}u\|_{L^2}^a\|u\|_{L^2}^{1-a}
\end{equation}
with $a=\frac{2r+3}{2r+4}<1$. Thus, if $k\geq 4$, by the last inequality and (\ref{i1}),
\begin{equation}\label{i1new}
I_1\leq C\sum_{r=2}^{[(k-2)/2]}\|D^{r+2}u\|_{L^2}^{2a}\|D^{k-r}u\|_{L^2}^2\|u\|_{L^2}^{2(1-a)}.
\end{equation}
Note that when $k\geq 4$, $r+2\leq k-1$ for all $r\leq [(k-2)/2]$. Hence, by induction hypothesis (\ref{ex:hyp}), it follows from (\ref{i1new})
\begin{equation}\label{i1final}
I_1\leq C_k \Psi^2_k(t).
\end{equation}
Combining (\ref{ineq:i}), (\ref{i1}), (\ref{ineq:imid}) and (\ref{i1final}) gives
\begin{equation}\label{ifinal}
I\leq C_k\Psi^2_k(t)
\end{equation}

For $II$, we have
\begin{align}\label{ineq:ii}
II&\leq 2\|\nabla d\|_{L^\infty}^2\|D^{k+1}d\|_{L^2}^2+2\|D^2 d\|_{L^\infty}^2\|D^{k}d\|_{L^2}^2+II_1,
\end{align}
with
\begin{equation}\label{ii1}
II_1=
\begin{cases}
&0, \ \ \mbox { if } k=2,3\\
&\sum_{r=2}^{[(k-2)/2]}\|D^{r+1}d\|_{L^\infty}^2\|D^{k-r+1}d\|_{L^2}^2, \ \ \mbox { if } k\geq 4.
\end{cases}
\end{equation}
Since $d$ is regular in the sense of (\ref{reg}), $\|\nabla d\|_{L^\infty}$ and $\|D^2 d\|_{L^\infty}$ are bounded. By induction hypothesis (\ref{ex:hyp}), $\|D^{k}d\|_{L^2}^2$ is bounded by $\Psi_k^2$. Thus,
\begin{equation}\label{ineq:iimid}
\|\nabla d\|_{L^\infty}^2\|D^{k+1}d\|_{L^2}^2+\|D^2 d\|_{L^\infty}^2\|D^{k}d\|_{L^2}^2\leq C\Psi^2_k(t).
\end{equation}
Applying Gagliardo-Nirenberg interpolation inequality,
\begin{equation}\notag
\|D^{r+1}d\|_{L^\infty}\leq C\|D^{r+3}d\|_{L^2}^a\|\nabla d\|_{L^2}^{2(1-a)}
\end{equation}
with $a=\frac{2r+3}{2r+4}<1$. Thus, if $k\geq 4$, by the last inequality and (\ref{ii1}),
\begin{equation}\label{ii1new}
II_1\leq C\sum_{r=2}^{[(k-2)/2]}\|D^{r+3}d\|_{L^2}^{2a}\|D^{k-r+1}d\|_{L^2}^2\|\nabla d\|_{L^2}^{2(1-a)}.
\end{equation}
Note that when $k\geq 4$, $r+3\leq k$ for all $r\leq [(k-2)/2]$. Hence, by induction hypothesis (\ref{ex:hyp}), it follows from (\ref{ii1new})
\begin{equation}\label{ii1final}
II_1\leq C_k\Psi^2_k(t)
\end{equation}
Combining (\ref{ineq:ii}), (\ref{ii1}), (\ref{ineq:iimid}) and (\ref{ii1final}) gives
\begin{equation}\label{iifinal}
II\leq C_k\Psi^2_k(t).
\end{equation}

For $III$, we have
\begin{align}\label{ineq:iii}
III&\leq \|u\|_{L^\infty}^2\|D^{k+1}d\|_{L^2}^2+\|\nabla u\|_{L^\infty}^2\|D^{k}d\|_{L^2}^2+\|D^2 u\|_{L^\infty}^2\|D^{k-1}d\|_{L^2}^2\\
&+\|\nabla d\|_{L^\infty}^2\|D^{k}u\|_{L^2}^2+\|D^2 d\|_{L^\infty}^2\|D^{k-1}u\|_{L^2}^2+III_1,\notag
\end{align}
with
\begin{equation}\label{iii1}
III_1=
\begin{cases}
&0, \ \ \mbox { if } k=2,3,4,\\
&\sum_{r=3}^{k-2}\|D^rd\|_{L^\infty}^2\|D^{k-r+1}u\|_{L^2}^2, \ \ \mbox { if } k\geq 5.
\end{cases}
\end{equation}
Recall $u$ and $d$ are regular by (\ref{reg}), $\|D^r u\|_{L^\infty}$ and $\|D^r d\|_{L^\infty}$ are bounded for $r=0,1,2$. By induction hypothesis (\ref{ex:hyp}), $\|D^{k-1}u\|_{L^2}^2$, $\|D^{k-1}d\|_{L^2}^2$ and $\|D^{k}d\|_{L^2}^2$ are bounded by $\Psi_k^2$. Thus, by (\ref{ineq:iii}),
\begin{equation}\label{ineq:iiimid}
III\leq C\Psi^2_k(t)+III_1.
\end{equation}
Applying Gagliardo-Nirenberg interpolation inequality,
\begin{equation}\label{inter:dr}
\|D^rd\|_{L^\infty}\leq C\|D^{r+2}d\|_{L^2}^a\|d-w_0\|_{L^2}^{1-a}
\end{equation}
with $a=\frac{2r+3}{2r+4}<1$. Thus, by (\ref{iii1})
\begin{equation}\label{iii1new}
III_1\leq C\sum_{r=3}^{k-2}\|D^{r+2}d\|_{L^2}^{2a}\|D^{k-r+1}u\|_{L^2}^2\|d-w_0\|_{L^2}^{2(1-a)}.
\end{equation}
Note that $r+2\leq k$ and $k-r+1\leq k-2$ for all $r\in[3,k-2]$, when $k\geq 5$. Hence, by induction hypothesis (\ref{ex:hyp}), it follows from (\ref{iii1new})
\begin{equation}\label{iii1final}
III_1\leq C_k\Psi_k^2.
\end{equation}
Combining 
(\ref{ineq:iiimid}) and (\ref{iii1final}) gives
\begin{equation}\label{iiifinal}
III\leq C\Phi^2_k(t)+C_k.
\end{equation}

Note that $|d|\leq 1$. For $IV$, we have
\begin{align}\label{ineq:iv}
IV&\leq C\|D^{k}d\|_{L^2}^2+\|\nabla d\|_{L^\infty}^2\|D^{k-1}d\|_{L^2}^2+\|D^2 d\|_{L^\infty}^2\|D^{k-2}d\|_{L^2}^2\\
&+\sum_{r=3}^{k-3}\|D^r d\|_{L^\infty}^2\|D^{k-r}d\|_{L^2}^2.\notag
\end{align}
By (\ref{inter:dr}) and induction hypothesis (\ref{ex:hyp}), it follows from (\ref{ineq:iv})
\begin{equation}\label{ivfinal}
IV\leq C_k \Psi_k^2(t).
\end{equation}

Combining (\ref{eq:dtphi1}), (\ref{ifinal}), (\ref{iifinal}), (\ref{iiifinal}) and (\ref{ivfinal}) yields
\begin{equation}\notag
\frac{d}{dt}\Phi_k^2(t)+\|D^{k+1} u\|_{L^2}^2+\|D^{k+2} d\|_{L^2}^2\leq C_k\Psi_k^2(t).
\end{equation}

 Summing  the last inequalities  for $k= 0,\dots,j$ yields 

\begin{equation}\label{1}
\frac{d}{dt}\Psi_j^2(t)+ \Psi_{j+1}^2(t)\leq C_j\Psi_j^2(t).
\end{equation}
Integrating (\ref{1})
yields by the inductive Hypothesis
\begin{equation} \label{2}
\Psi_j^2(t) + \int_0^t \Psi_{j+1}^2(t)\leq \Psi_j^2(0)+ C_j \int_0^t \Psi_j^2(t)\leq C_j \left( \|u_0\|_{H^j(\R^3)} +\|d_0-\omega_0\|_{H^{j+1}(\R^3)} \right).
\end{equation}
This concludes the proof of the Theorem.
\cbdu

\bigskip

\section{Improvement of the velocity  decay rate}\label{sec:3}

In \cite{DQS}, for solutions with $u_0\in L^2(\R^3)\cap L^1(\R^3)$ , appropriate assumptions on the data $d_0$   and sufficiently small data, we showed that  the velocity satisfies $\|u\|_{L^2}^2\leq C(1+t)^{-\frac{1}{2}}$. In this section,
with the same assumptions  we show that the velocity has the optimal decay rate
$\|u\|_{L^2}^2\leq C(1+t)^{-\frac{3}{2}}$. Namely, we show

\begin{Lemma}\label{le:udecay} Assume $u_0\in H^1(\R^3)\cap L^1(\R^3)$ and  $d_0- w_0 \in H^{2}(\R^3)\cap L^1(\R^3)$. Let $u$ be a smooth solution obtained in Theorem \ref{exist} with data satisfying (\ref{smalldataH2}).  Then, $u$ satisfies 
\begin{equation}\label{eq:udecay}
\|u\|_{L^2(\R^3)}^2\leq C(1+t)^{-\frac{3}{2}}.
\end{equation}
\end{Lemma}
\pf
As in \cite{DQS}, the Fourier splitting method, \cite{Sch1,Sch2}, is applied to establish the decay of solution. We decompose the frequency domain $\R^3$ into two time dependent subdomains $S(t)$ and its complement $S^c(t)$. Define 
\begin{equation}\label{S}
S(t)=\left\{\xi\in\R^3:|\xi|\leq r(t)=(\frac{k}{1+t})^{1/2}\right\}
\end{equation}
for a constant $k$ that will be specified bellow. One of the key estimates to establish decay of velocity 
\begin{equation}\label{decay:u1}
\|u\|_{L^2(\R^3)}^2\leq C(1+t)^{-\frac{1}{2}}
\end{equation}
in \cite{DQS} is that the Fourier transform of the solution $u$ satisfies
\begin{equation}\notag
|\hat u(\xi,t)|\leq C|\xi|^{-1} \ \ \ \mbox{ for } \xi\in S(t).
\end{equation}
 We show first  that with the decay estimate (\ref{decay:u1}), the Fourier transform of
the velocity $u$ satisfies
\begin{equation}\label{est:fourieru}
|\hat u(\xi,t)|\leq C, \ \ \ \mbox{ for } \xi\in S(t)
\end{equation}
for an absolute constant $C$. As a consequence, we proceed with a similar analysis as in \cite{DQS} and obtain the optimal decay (\ref{eq:udecay}) for the velocity.\\
Taking the Fourier transform of Navier-Stokes equation in system (\ref{LCD}) yields
\bg\label{FNSE}
\hat u_t+|\xi|^2\hat u=G(\xi,t)
\ed
where
\bg\notag
G(\xi,t)=-\mathcal F(u\cdot\nabla u)-\mathcal F(\nabla \pi)-\mathcal F(\nabla\cdot(\nabla d\otimes\nabla d)),
\ed
and $\mathcal F$ indicates the Fourier transform. Multiplying (\ref{FNSE}) by the integrating factor $e^{|\xi|^2t}$ yields
\bg\notag
\frac{d}{dt}[e^{|\xi|^2t}\hat u]=e^{|\xi|^2t}G(\xi,t).
\ed
Integrating in time gives
\bg\label{Fueq}
\hat u(\xi,t)=e^{-|\xi|^2t}\hat u_0+\int_0^te^{-|\xi|^2(t-s)}G(\xi,s)ds.
\ed
We analyze each term in $G(\xi,t)$ separately. We have
\bg\label{Ftensoru}
|\mathcal F(u\cdot\nabla u)|=|\mathcal F(\nabla\cdot(u\otimes u))|\leq \sum_{i,j}\int_{\R^3}|u^iu^j||\xi_j|dx\leq C(1+t)^{-1/2}|\xi|
\ed
due to (\ref{decay:u1}).\\
In \cite{DQS}, we obtained decay for the gradient of director field $d$ as
\begin{equation}\notag
\|\nabla d\|_{L^2(\R^3)}^2\leq C(1+t)^{-3/4}.
\end{equation}
Thus, we have
\bg\label{Ftensord}
|\mathcal F(\nabla\cdot(\nabla d\otimes\nabla d))|\leq C(1+t)^{-3/4}|\xi|.
\ed
Taking divergence of the velocity equation in system (\ref{LCD}) yields
\bg\notag
\Delta \pi=-\sum_{i,j}\frac{\partial^2}{\partial x_i\partial x_j}(u^iu^j)-\sum_{i,j}\frac{\partial^2}{\partial x_i\partial x_j}(\nabla d^i\nabla d^j).
\ed
Taking the Fourier transform of the last equation gives
\bg\notag
|\xi|^2\mathcal F(\pi)=-\sum_{i,j}\xi_i\xi_j\mathcal F(u^iu^j)-\sum_{i,j}\xi_i\xi_j\mathcal F(\nabla d^i\nabla d^j).
\ed
Combining (\ref{Ftensoru}), (\ref{Ftensord}) and the last equation yields
\bg\notag
\mathcal F(\pi)\leq C(1+t)^{-1/2},
\ed
and thus 
\begin{equation}\label{Fp}
\mathcal F(\nabla \pi)\leq C(1+t)^{-1/2}|\xi|.
\end{equation} 
Combining (\ref{Ftensoru}), (\ref{Ftensord}) and (\ref{Fp}) yields
\bg\label{G}
|G(\xi,t)|\leq C(1+t)^{-1/2}|\xi|, \ \ \mbox{ for } \xi\in S(t).
\ed
From (\ref{Fueq}) and (\ref{G}) we have
\bg\label{Fuineq}
|\hat u(\xi,t)|\leq e^{-|\xi|^2t}|\hat u_0|+\int_0^te^{-|\xi|^2(t-s)}(1+s)^{-1/2}|\xi|ds.
\ed
Since $u_0\in L^1$, we have $|\hat u_0|\leq C$ for all $\xi$. Performing integration on the right hand side of (\ref{Fuineq}) gives
\bg\notag
|\hat u(\xi,t)|\leq Ce^{-|\xi|^2t}+C(1+t)^{1/2}|\xi|\leq C
\ed
since $\xi\in S(t)$ and $|\xi|\leq C(1+t)^{-1/2}$. This completes the proof of (\ref{est:fourieru}).\\

Multiplying the velocity equation in (\ref{LCD}) by $u$ and integrating by parts yields
\begin{align}\notag
\frac{1}{2}\frac{d}{dt}\int_{\R^3}|u|^2dx+\int_{\R^3}|\nabla u|^2dx&=\int_{\R^3}\nabla u(\nabla d\otimes\nabla d)dx\\
&\leq \frac{1}{2}\int_{\R^3}|\nabla u|^2dx+C\int_{\R^3}|\nabla d\otimes\nabla d|^2dx.\notag
\end{align}

Thus, 
\bg\label{NSEen}
\frac{d}{dt}\int_{\R^3}|u|^2+\int_{\R^3}|\nabla u|^2dx\leq C\int_{\R^3}|\nabla d\otimes\nabla d|^2dx.
\ed
The right hand side of (\ref{NSEen}) is estimated by
\begin{align}\notag
\int_{\R^3}|\nabla d\otimes\nabla d|^2dx&=\int_{\R^3}(\nabla d\otimes\nabla d)(\nabla d\otimes\nabla d)dx\\
&=-3\int_{\R^3}(d-w_0)\;\Delta d\;\nabla d\otimes\nabla d\notag\\
&\leq \frac{1}{2}\int_{\R^3}|\nabla d\otimes\nabla d|^2dx+C\int_{\R^3}|d-w_0|^2|\Delta d|^2dx\notag.
\end{align}
It follows that
\begin{align}\notag
\int_{\R^3}|\nabla d\otimes\nabla d|^2dx&\leq C\int_{\R^3}|d-w_0|^2|\Delta d|^2dx\\
&\leq C\left(\int_{\R^3}|d-w_0|^pdx\right)^{\frac{2}{p}}\left(\int_{\R^3}|\Delta d|^{\frac{2p}{p-2}}dx\right)^{\frac{p-2}{p}}\notag\\
&=C\|d-w_0\|_{L^p(\R^3)}^2\left(\int_{\R^3}|\Delta d|^{2+\frac{4}{p-2}}dx\right)^{\frac{p-2}{p}}\notag\\
&\leq C\|d-w_0\|_{L^p(\R^3)}^2\notag,
\end{align}
for $p> 2$. The last step followed from the energy estimate (\ref{ex:lady}) and recalling that $\|\Delta d\|_{L^\infty(\R^3\times[0,T])}$ is bounded since $d\in H^m(\R^3)$ is regular by (\ref{reg}) in Theorem \ref{exist}. Thus, (\ref{ddecay}) in Theorem \ref{Mthm} yields
\bg\notag
\int_{\R^3}|\nabla d\otimes\nabla d|^2dx\leq C(1+t)^{-3(1-\frac{1}{p})},
\ed
for any $p\geq 2$. Therefore, it follows from (\ref{NSEen})
\bg\label{NSEen1}
\frac{d}{dt}\int_{\R^3}|u|^2dx+\int_{\R^3}|\nabla u|^2dx\leq C(1+t)^{-3(1-\frac{1}{p})}.
\ed

Applying Plancherel's theorem to (\ref{NSEen1}) gives
\bg\notag
\frac{d}{dt}\int_{\R^3}|\hat{u}|^2d\xi+\int_{\R^3}|\xi|^2|\hat{u}|^2d\xi\leq C(1+t)^{-3(1-\frac{1}{p})}.
\ed
We reorganize the last inequality  as
\begin{align}\label{NSEF}
\frac{d}{dt}\int_{\R^3}|\hat{u}|^2d\xi&\leq-\int_{S(t)^c}|\xi|^2|\hat{u}|^2d\xi-\int_{S(t)}|\xi|^2|\hat{u}|^2d\xi
+C(1+t)^{-3(1-\frac{1}{p})}\\
&\leq-\frac{k}{1+t}\int_{S(t)^c}|\hat{u}|^2d\xi-\int_{S(t)}|\xi|^2|\hat{u}|^2d\xi
+C(1+t)^{-3(1-\frac{1}{p})}\notag\\
&\leq-\frac{k}{1+t}\int_{\R^3}|\hat{u}|^2d\xi+
\frac{k}{1+t}\int_{S(t)}|\hat{u}|^2d\xi+C(1+t)^{-3(1-\frac{1}{p})}\notag.
\end{align}
Due to (\ref{est:fourieru}) we have $|\hat u|\leq C$, hence
\bg\notag
\int_{S(t)}|\hat u|^2d\xi\leq C\int_0^{r(t)}r^2dr\leq C(1+t)^{-3/2}.
\ed
By  (\ref{NSEF})  choosing any $p\geq6$
\begin{align}\notag
\frac{d}{dt}\int_{\R^3}|\hat{u}|^2d\xi+\frac{k}{1+t}\int_{\R^3}|\hat{u}|^2d\xi&\leq
C(1+t)^{-\frac{5}{2}}+C(1+t)^{-3(1-\frac{1}{p})}\\
&\leq C(1+t)^{-\frac{5}{2}}\notag.
\end{align}

Multiplycation  by the integrating factor $(1+t)^k$ yields
\bg\notag
\frac{d}{dt}\left[(1+t)^k\int_{\R^3}|\hat{u}|^2d\xi\right]\leq C(1+t)^{k-\frac{5}{2}}.
\ed
Integration  in time  gives
\bg\notag
(1+t)^{k}\int_{\R^3}|\hat{u}|^2d\xi\leq\int_{\R^3}|\hat{u}(\xi,0)|^2d\xi+C[(1+t)^{k-\frac{3}{2}}-1].
\ed
Thus,
\bg\notag
\int_{\R^3}|\hat{u}|^2d\xi\leq(1+t)^{-k}\int_{\R^3}|\hat{u}(\xi,0)|^2d\xi+C[(1+t)^{-\frac{3}{2}}-(1+t)^{-k}].
\ed
Since $u_0\in L^2$, $\hat u(0)\in L^2$ by Plancherel's theorem. We choose $k>3/2$.
Hence
\bg\notag
\int_{\R^3}|u|^2dx\leq C(1+t)^{-\frac{3}{2}}.
\ed
The proof of the Lemma is complete.
\cbdu

The rest of this section deals with an auxiliary estimate for  the velocity $u$ in $L^\infty$.

\begin{Lemma}\label{le:u-infty-decay}
Let $u$ be the solution of system (\ref{LCD}) with initial data satisfying the conditions in
Theorem \ref{exist}. Then we have
\begin{equation}\label{decay:uinfty}
\|u\|_{L^\infty(\R^3)}\leq C_m(1+t)^{-\frac{3}{4}(1-\frac{3}{2m})},
\end{equation}
for $m\geq 1$.
\end{Lemma}
\pf
By Gagliardo-Nirenberg inequality (\ref{GN}), we have
\begin{equation}\notag
\|u\|_{L^\infty}\leq \|D^mu\|_{L^2}^a\|u\|_{L^2}^{1-a}
\end{equation}
with $a=\frac{3}{2m}$. Combining the estimates (\ref{ex:Hlady}),  (\ref{eq:udecay}) and the 
above inequality yields
\begin{equation}\notag
\|u\|_{L^\infty}\leq C_m(1+t)^{-\frac{3}{4}(1-\frac{3}{2m})}
\end{equation}
which establishes the conclusion of the lemma.
\cbdu

\begin{Remark}
The rate in Lemma \ref{le:u-infty-decay} will be improved to $(1+t)^{-3/2}$ in Section \ref{sec:5}.
\end{Remark}

\bigskip

\section{Optimal decay rate of $\nabla d$}\label{sec:4}

The goal of this section is to  establish that $\nabla d$ decays in $L^2(\R^3)$ at the rate $(1+t)^{-\frac{5}{4}}$.\\

We first obtain an auxiliary decay rate for $d-w_0$ in $L^\infty(\R^3)$,

\begin{Lemma}\label{le:d-infty}
Let $d$ be the solution obtained in Theorem \ref{exist}. Assume the initial data satisfies the conditions in Theorem \ref{Mthm}. Then 
\bg\label{ineq:d-infty-decay}
\|d(\cdot,t)-w_0\|_{L^\infty(\R^3)}\leq C(1+t)^{-\frac{3}{2}\cdot\frac{p-1}{p+2}},  \ \ p>1
\ed
where the constant $C$ depends only on the initial data.
\end{Lemma}
\pf
By the Gagliardo-Nirenberg interpolation inequality,
\bg\notag
\|d-w_0\|_{L^\infty}\leq C\|D^3 d\|_{L^2}^a\|d-w_0\|_{L^p}^{1-a}
\ed
with $a=\frac{2}{2+p}$. Due to the energy estimate (\ref{ex:lady}) and the decay estimate in 
Theorem \ref{Mthm}, we have that
\bg\notag
\|d-w_0\|_{L^\infty}\leq C(1+t)^{-\frac{3}{2}(1-\frac{1}{p})(1-\frac{2}{2+p})},
\ \ p>1,
\ed
which proves (\ref{ineq:d-infty-decay}).
\cbdu

\begin{Corollary}\label{cor1}
Let $d$ be the solution obtained in Theorem \ref{exist}. Assume the initial data satisfies the conditions in Theorem \ref{Mthm}. Then 
\bg\label{ineq:dclose}
(d+w_0)\cdot d\geq 0,  
\ed
for sufficiently large time $t$.
\end{Corollary}
\pf
Taking $p=7$ in (\ref{ineq:d-infty-decay}) yields 
$$
\|d-w_0\|_{L^\infty}\leq C(1+t)^{-1}.
$$
Hence, we have 
\begin{align}\notag
(d-w_0)\cdot d&=|d|^2+d\cdot w_0\\
&=|d|^2+(d-w_0)\cdot w_0+1\notag\\
&\geq |d|^2+1-|d-w_0||w_0|\notag\\
&\geq |d|^2+1-C(1+t)^{-1}\notag\\
&\geq 0\notag
\end{align}
for sufficiently large time $t$. This proves (\ref{ineq:dclose}).
\cbdu

\begin{Lemma}\label{le:gradient-d-decay}
Let $d$ be the solution obtained in Theorem \ref{exist}. Assume the initial data satisfies the conditions in Theorem \ref{Mthm}. Then 
\bg\label{ineq:d-w0decay}
\|\nabla d(\cdot,t)\|_{L^2(\R^3)}\leq C(1+t)^{-\frac{5}{4}}.
\ed
The constant $C$ depends only on the initial data.
\end{Lemma}
\pf
Multiply the director field equation in (\ref{LCD}) by $\Delta d$ and integrate over $\R^3$, then
\bg\label{ineq:grad-d}
\frac{1}{2}\frac{d}{dt}\int_{\R^3}|\nabla d|^2dx+\int_{\R^3}|\Delta d|^2dx=\int_{\R^3}(u\cdot\nabla d)\Delta d dx+\int_{\R^3}f(d)\Delta d dx.
\ed
We need the following auxiliary estimates
\begin{align}\notag
\left |\int_{\R^3}(u\cdot\nabla d)\Delta d dx\right |&\leq\frac{1}{4}\int_{\R^3}|\Delta d|^2dx+C\|u\|_{L^\infty}^2\int_{\R^3}|\nabla d|^2dx
\end{align}
and since $|w_0|=1$
\begin{align}\notag
\int_{\R^3}f(d)\Delta d dx &=-\frac{1}{\eta^2}\int_{\R^3}\nabla[(d+w_0)(d-w_0)d]\nabla ddx\\
&\leq C\|d-w_0\|_{L^\infty}\int_{\R^3}|\nabla d|^2dx-\frac{1}{\eta^2}\int_{\R^3}(d+w_0)\cdot d|\nabla d|^2dx\notag\\
&\leq C\|d-w_0\|_{L^\infty}\int_{\R^3}|\nabla d|^2dx\notag,
\end{align}
where we used the facts that $|d|\leq 1$ and $(d+w_0)\cdot d\geq 0$ for large time, see Corollary \ref{cor1}. Combining the last two inequalities with (\ref{ineq:grad-d}) yields
\begin{align}\notag
&\frac{d}{dt}\int_{\R^3}|\nabla d|^2dx+\int_{\R^3}|\Delta d|^2dx\\
&\leq C(\|u\|_{L^\infty}^2+\|d-w_0\|_{L^\infty})\int_{\R^3}|\nabla d|^2dx\notag\\
&\leq C((1+t)^{-\frac{3}{2}(1-\frac{3}{2m})}+(1+t)^{-\frac{3}{2}\cdot\frac{p-1}{p+2}})\int_{\R^3}|\nabla d|^2dx\notag
\end{align}
due to Lemma \ref{le:u-infty-decay} and Lemma \ref{le:d-infty}. Taking $m=5, p=7$ yields that
\begin{equation}\label{ineq:grad-d-energy}
\frac{d}{dt}\int_{\R^3}|\nabla d|^2dx+\int_{\R^3}|\Delta d|^2dx\leq C_0(1+t)^{-1}\int_{\R^3}|\nabla d|^2dx.
\end{equation}
Apply now the Fourier splitting method to obtain the decay estimate. We choose an appropriate constant $k$ for the time dependent sphere $S(t)$ as in (\ref{S}), such that $k-C_0>5/2$. Here the constant $C_0$ is the one on the right hand side of (\ref{ineq:grad-d-energy}). Proceeding by the Fourier splitting method gives 
\begin{align}\label{ineq:delta-d-s}
\int_{\R^3}|\Delta d|^2dx&\geq \int_{\R^3\setminus S}|\xi|^2|\mathcal F(\nabla d)|^2d\xi\\
&\geq \frac{k}{1+t}\int_{\R^3\setminus S}|\mathcal F(\nabla d)|^2d\xi\notag\\
&=\frac{k}{1+t}\int_{\R^3}|\mathcal F(\nabla d)|^2d\xi- \frac{k}{1+t}\int_{S}|\mathcal F(\nabla d)|^2d\xi\notag.
\end{align}
Combining (\ref{ineq:grad-d-energy}) and (\ref{ineq:delta-d-s}) yields, 
\begin{align}\label{ineq:grad-d-decay}
&\frac{d}{dt}\int_{\R^3}|\nabla d|^2dx+\frac{k-C_0}{1+t}\int_{\R^3}|\nabla d|^2dx\\
&\leq \frac{k}{1+t}\int_{S}|\mathcal F(\nabla d)|^2d\xi\notag\\
&\leq \frac{C}{(1+t)^2}\int_{S}|\mathcal F(d-w_0)|^2d\xi\notag.
\end{align}
Multiplying (\ref{ineq:grad-d-decay}) by the factor $(1+t)^{k-C_0}$ and integrating in time yields 
\begin{align}\notag
\int_{\R^3}|\nabla d|^2dx&\leq (1+t)^{-(k-C_0)}\int_{\R^3}|\nabla d_0|^2dx+(1+t)^{-1}\int_{\R^3}|d-w_0|^2dx\\
&\leq C(1+t)^{-\frac{5}{2}}\notag,
\end{align}
due to the decay estimate in Theorem \ref{Mthm} and the choice $k-C_0>5/2$. This completes the proof of the Lemma.
\cbdu

\bigskip

\section {Decay of solutions in higher order Sobolev spaces}
\label{sec:5}

In this section we obtain the decay estimates for $u$ in $H^m(\R^3)$ and $d$ in $H^{m}(\R^3)$ with $m\geq 1$. The method involves induction and a Fourier
splitting argument.\\

We recall  the following Gagliardo-Nirenberg inequalities
\bg\label{ineq:interp-u}
\|D^iu\|_{L^\infty}\leq C\|D^{m+1}u\|_{L^2}^{a_i}\|u\|_{L^2}^{1-a_i}
\ed
\bg\label{ineq:interp-d}
\|D^id\|_{L^\infty}\leq C\|D^{m+1}d\|_{L^2}^{a_i}\|d-w_0\|_{L^2}^{1-a_i}
\ed
with $a_i=\frac{i+3/2}{m+1}$, for $i\geq 1$ and $1+3/2<m-1/2$. .\\

We first establish the following auxiliary estimates.
\begin{Lemma}\label{le:infty-interp}
Let $d$ be the solution obtained in Theorem \ref{exist}. Then $d$ satisfies 
\bg\label{dec:d-infty}
\|\nabla d\|_{L^\infty}\leq C(1+t)^{-\frac{5}{4}(1-\frac{3}{2k})}
\ed
\bg\label{dec:d2-infty}
\|D^2 d\|_{L^\infty}\leq C(1+t)^{-\frac{5}{4}(1-\frac{5}{2k})}
\ed
for $k\geq 1$.
\end{Lemma}
\pf
Taking $w=\nabla d$, $k=0$, $r=\infty$, and $p=q=2$ in (\ref{GN}) yields, for $k\geq 2$
\begin{align}\notag
\|\nabla d\|_{L^\infty}&\leq C\|D^{k}d\|_{L^2}^{a}\|\nabla d\|_{L^2}^{1-a}\\
&\leq C(1+t)^{-\frac{5}{4}(1-\frac{3}{2k})}\notag
\end{align}
due to the higher order energy estimate (\ref{ex:lady}) in Theorem \ref{exist} and the decay estimate (\ref{ineq:d-w0decay}) in Lemma \ref{le:gradient-d-decay}. The constant $C$ depends only on initial data. Thus (\ref{dec:d-infty}) is proved. The proof of (\ref{dec:d2-infty}) is similar and as such is omitted.
\cbdu
\begin{Remark}
The decay of $\nabla d$ and $\Delta d$ in $L^\infty(\R^3)$ will be improved to the optimal rates.
\end{Remark}

We now establish a higher order energy estimate for the solution. The ideas are based on work in \cite{Sch94}.

\begin{Lemma}\label{le:hm-estimate}
Let $(u,d)$ be the solution obtained in Theorem \ref{exist}. Assume the initial data satisfies the conditions in Theorem \ref{Mthm}. Then $(u,d)$ satisfies, for $m\geq 1$
\begin{align}\label{ineq:hm-estimate}
&\frac{d}{dt}\int_{\R^3}|D^mu|^2+|D^{m}d|^2dx+\int_{\R^3}|D^{m+1}u|^2+|D^{m+1}d|^2dx\\
&\leq C_m(\|u\|_{L^\infty}^2+\|d-w_0\|_{L^\infty}^2+\|D^2 d\|_{L^\infty}^2)\int_{\R^3}|D^mu|^2+|D^{m}d|^2dx\notag\\
&+C_m\|d-w_0\|_{L^\infty}^2\|D^{m-1}d\|_{L^2}^2+R_m\notag,
\end{align}
where
\begin{equation}\label{R}
R_m=
\begin{cases}
0, \ \ \ m=1,2;\\
\sum_{i=1}^{\frac{m}{2}}\|u\|_{L^2}^{2}\|D^{m-i}u\|_{L^2}^{\frac{2}{1-a_i}}+
\sum_{i=2}^{\frac{m-2}{2}}\|d-w_0\|_{L^2}^{2}\|D^{m-i+1}d\|_{L^2}^{\frac{2}{1-a_{i+1}}}\\
+\sum_{i=1}^{m-1}\|d-w_0\|_{L^2}^2\|D^{m-i}u\|_{L^2}^{\frac{2}{1-a_{i}}}+
\sum_{i=1}^{\frac{m-1}{2}}\|d-w_0\|_{L^2}^{2}\|D^{m-i-1}d\|_{L^2}^{\frac{2}{1-a_i}},  \\
\equiv R_{m1}+R_{m2}+R_{m3}+R_{m4}, 
\ \ \mbox { for } m\geq 3
\end{cases}
\end{equation}
with $a_i=\frac{i+3/2}{m+1}$.
The constants $C$ in the inequality depend only on the initial data.
\end{Lemma}
\pf
Taking the $m$-th derivative on the first equation in (\ref{LCD}),
 multiplying it by $D^mu$ and integrating over $\R^3$ yields
\begin{align}\notag
&\frac{1}{2}\frac{d}{dt}\int_{\R^3}|D^mu|^2dx+\int_{\R^3}|D^{m+1}u|^2dx\\
&=-\int_{\R^3}D^m(u\cdot\nabla u)D^mudx-\int_{\R^3}D^m(\nabla\cdot(\nabla d\otimes\nabla d))D^mudx\notag\\
&=\int_{\R^3}D^{m-1}(u\cdot\nabla u)D^{m+1}udx+\int_{\R^3}D^m(\nabla d\otimes\nabla d)D^{m+1}udx\notag\\
&\leq\frac{1}{4}\int_{\R^3}|D^{m+1}u|^2dx+\int_{\R^3}|D^{m-1}(u\cdot\nabla u)|^2dx+
\int_{\R^3}|D^m(\nabla d\otimes\nabla d)|^2dx\notag.
\end{align}
Therefore, 
\begin{align}
\label{ineq:mth-energy}
&\frac{d}{dt}\int_{\R^3}|D^mu|^2dx+\frac{3}{2}\int_{\R^3}|D^{m+1}u|^2dx\\
&\leq\int_{\R^3}|D^{m-1}(u\cdot\nabla u)|^2dx+
\int_{\R^3}|D^m(\nabla d\otimes\nabla d)|^2dx\notag\\
&\leq\|u\|_{L^\infty}^2\|D^mu\|_{L^2}^2+\|\nabla d\|_{L^\infty}^2\|D^{m+1}d\|_{L^2}^2+\|D^2 d\|_{L^\infty}^2\|D^{m}d\|_{L^2}^2\notag\\
&+\sum_{1\leq i\leq m/2}\|D^iu\|_{L^\infty}^2\|D^{m-i}u\|_{L^2}^2
+\sum_{2\leq i\leq (m-2)/2}\|D^{i+1}d\|_{L^\infty}^2\|D^{m-i+1}d\|_{L^2}^2\notag.
\end{align}
Using (\ref{ineq:interp-u}) gives that
\begin{align}\label{ineq:ui}
&\sum_{1\leq i\leq m/2}\|D^iu\|_{L^\infty}^2\|D^{m-i}u\|_{L^2}^2\\
&\leq C\sum_{1\leq i\leq m/2}\|D^{m+1}u\|_{L^2}^{2a_i}\|u\|_{L^2}^{2(1-a_i)}\|D^{m-i}u\|_{L^2}^2\notag\\
&\leq \frac{1}{8}\|D^{m+1}u\|_{L^2}^2+C\sum_{1\leq i\leq m/2}\|u\|_{L^2}^{2}\|D^{m-i}u\|_{L^2}^{2/(1-a_i)}\notag.
\end{align}
Using (\ref{ineq:interp-d}) gives that
\begin{align}\label{ineq:di}
&\sum_{2\leq i\leq (m-2)/2}\|D^{i+1}d\|_{L^\infty}^2\|D^{m-i+1}d\|_{L^2}^2\\
&\leq C\sum_{2\leq i\leq (m-2)/2}\|D^{m+1}d\|_{L^2}^{2a_{i+1}}\|d-w_0\|_{L^2}^{2(1-a_{i+1})}\|D^{m-i+1}d\|_{L^2}^2\notag\\
&\leq \frac{1}{8}\|D^{m+1}d\|_{L^2}^2+C\sum_{2\leq i\leq (m-2)/2}\|d-w_0\|_{L^2}^{2}\|D^{m-i+1}d\|_{L^2}^{2/(1-a_{i+1})}\notag.
\end{align}

Taking the $m$-th derivative on the second equation in (\ref{LCD}),
 multiplying it by $D^{m}d$ and integrating over $\R^3$ yields
\begin{align}\notag
&\frac{1}{2}\frac{d}{dt}\int_{\R^3}|D^{m}d|^2dx+\int_{\R^3}|D^{m+1}d|^2dx\\
&=-\int_{\R^3}D^{m}(u\cdot\nabla d)D^{m}d+\frac{1}{\eta^2}D^{m}[(|d|^2-1)d]D^{m}ddx\notag\\
&=\int_{\R^3}D^{m}(u\otimes d)D^{m+1}d+\frac{1}{\eta^2}D^{m-1}[(|d|^2-1)d]D^{m+1}ddx\notag\\
&\leq\frac{1}{4}\int_{\R^3}|D^{m+1}d|^2dx+\int_{\R^3}|D^{m}(u\otimes d)|^2+C|D^{m-1}[(|d|^2-1)d]|^2dx\notag.
\end{align}
Therefore,
\begin{align}\label{ineq:m1th-energy}
&\frac{d}{dt}\int_{\R^3}|D^{m}d|^2dx+\frac{3}{2}\int_{\R^3}|D^{m+1}d|^2dx\\
&\leq C\int_{\R^3}|D^{m}(u\otimes d)|^2+|D^{m-1}[(|d|^2-1)d]|^2dx\notag.
\end{align}
The first integral on the right hand side of (\ref{ineq:m1th-energy}) is estimated as
\begin{align}\label{ineq:udi}
&\int_{\R^3}|D^{m}(u\otimes d)|^2dx\\
&\leq \|u\|_{L^\infty}^2\|D^{m}d\|_{L^2}^2+\|d-w_0\|_{L^\infty}^2\|D^{m}u\|_{L^2}^2\notag\\
&+\sum_{i=1}^{m-1}\|D^{i}d\|_{L^\infty}^2\|D^{m-i}u\|_{L^2}^2\notag\\
&\leq \|u\|_{L^\infty}^2\|D^{m}d\|_{L^2}^2+\|d-w_0\|_{L^\infty}^2\|D^{m}u\|_{L^2}^2\notag\\
&+\sum_{i=1}^{m-1}\|D^{m+1}d\|_{L^2}^{2a_{i}}\|d-w_0\|_{L^2}^{2(1-a_{i})}\|D^{m-i}u\|_{L^2}^2\notag\\
&\leq \|u\|_{L^\infty}^2\|D^{m}d\|_{L^2}^2+\|d-w_0\|_{L^\infty}^2\|D^{m}u\|_{L^2}^2\notag\\
&+\frac{1}{8}\|D^{m+1}d\|_{L^2}^2+C_m\sum_{i=1}^{m-1}\|d-w_0\|_{L^2}^2\|D^{m-i}u\|_{L^2}^{2/(1-a_{i})}\notag.
\end{align}
The second integral on the right hand side of (\ref{ineq:m1th-energy}) is estimated as
\begin{align}\label{ineq:fdm}
&\int_{\R^3}|D^{m-1}[(|d|^2-1)d]|^2dx\\
&=\int_{\R^3}\sum_{i=0}^{m-1}|D^i[(d+w_0)(d-w_0)]D^{m-i-1}d|^2dx\notag\\
&\leq\|d-w_0\|_{L^\infty}^2\|D^{m-1}d\|_{L^2}^2+(1+\|d-w_0\|_{L^\infty}^2)\sum_{i=1}^{(m-1)/2}\|D^id\|_{L^\infty}^2\|D^{m-i-1}d\|_{L^2}^2\notag\\
&+\sum_{i=1}^{(m-1)/2}\sum_{j=1}^{i-1}\|D^jd\|_{L^\infty}^2\|D^{i-j}d\|_{L^\infty}^2\|D^{m-i-1}d\|_{L^2}^{2}\notag.
\end{align}
Using the interpolation inequality (\ref{ineq:interp-d}) gives that
\begin{align}\notag
&\sum_{i=1}^{(m-1)/2}\|D^id\|_{L^\infty}^2\|D^{m-i-1}d\|_{L^2}^2\\
&\leq C\sum_{i=1}^{(m-1)/2}\|D^{m+1}d\|_{L^2}^{2a_i}\|d-w_0\|_{L^2}^{2(1-a_i)}\|D^{m-i-1}d\|_{L^2}^2\notag\\
&\leq \frac{1}{8}\|D^{m+1}d\|_{L^2}^2+C_m\sum_{i=1}^{(m-1)/2}\|d-w_0\|_{L^2}^{2}\|D^{m-i-1}d\|_{L^2}^{2/(1-a_i)}\notag.
\end{align}
Using the interpolation inequality (\ref{ineq:interp-d}) again to the last term in (\ref{ineq:fdm}) gives that
\begin{align}\notag
&\sum_{i=1}^{(m-1)/2}\sum_{j=1}^{i-1}\|D^jd\|_{L^\infty}^2\|D^{i-j}d\|_{L^\infty}^2\|D^{m-i-1}d\|_{L^2}^2\\
&\leq C\sum_{i=1}^{(m-1)/2}\sum_{j=1}^{i-1}\|D^{m+1}d\|_{L^2}^{2(a_j+a_{i-j})}\|d-w_0\|_{L^2}^{2(2-a_j-a_{i-j})}\|D^{m-i-1}d\|_{L^2}^2\notag\\
&\leq \frac{1}{8}\|D^{m+1}d\|_{L^2}^2+C_m\sum_{i=1}^{(m-1)/2}\|d-w_0\|_{L^2}^{\frac{2(2-a_i)}{1-a_i}}\|D^{m-i-1}d\|_{L^2}^{2/(1-a_i)}\notag.
\end{align}
The last two inequalities combined with (\ref{ineq:fdm}) imply that
\begin{align}\label{ineq:new-fdm}
&\int_{\R^3}|D^{m-1}[(|d|^2-1)d]|^2dx\\
&\leq\frac{1}{4}\|D^{m+1}d\|_{L^2}^2+\|d-w_0\|_{L^\infty}^2\|D^{m-1}d\|_{L^2}^2\notag\\
&+(1+\|d-w_0\|_{L^\infty}^2)\sum_{i=1}^{(m-1)/2}\|d-w_0\|_{L^2}^{2}\|D^{m-i-1}d\|_{L^2}^{2/(1-a_i)}\notag,
\end{align}
where we used the facts that $\frac{2-a_i}{1-a_i}>1$ and that $\|d-w_0\|_{L^2}$  decays for large time.\\

Finally, adding the two inequalities (\ref{ineq:mth-energy}) and (\ref{ineq:m1th-energy}), and using (\ref{ineq:ui}), (\ref{ineq:di}), (\ref{ineq:udi}) and (\ref{ineq:new-fdm}) gives that
\begin{align}
&\frac{d}{dt}\int_{\R^3}|D^mu|^2+|D^{m}d|^2dx+\int_{\R^3}|D^{m+1}u|^2+|D^{m+1}d|^2dx\notag\\
&\leq C(\|u\|_{L^\infty}^2+\|d-w_0\|_{L^\infty}^2+\|D^2 d\|_{L^\infty}^2)\int_{\R^3}|D^mu|^2+|D^{m}d|^2dx\notag\\
&+C_m\|d-w_0\|_{L^\infty}^2\|D^{m-1}d\|_{L^2}^2+R_m\notag,
\end{align}
where
\begin{equation}\notag
R_m=
\begin{cases}
0, \ \ \ m=1,2;\\
\sum_{i=1}^{\frac{m}{2}}\|u\|_{L^2}^{2}\|D^{m-i}u\|_{L^2}^{\frac{2}{1-a_i}}+
\sum_{i=2}^{\frac{m-2}{2}}\|d-w_0\|_{L^2}^{2}\|D^{m-i+1}d\|_{L^2}^{\frac{2}{1-a_{i+1}}}\\
+\sum_{i=1}^{m-1}\|d-w_0\|_{L^2}^2\|D^{m-i}u\|_{L^2}^{\frac{2}{1-a_{i}}}+
\sum_{i=1}^{\frac{m-1}{2}}\|d-w_0\|_{L^2}^{2}\|D^{m-i-1}d\|_{L^2}^{\frac{2}{1-a_i}},  \\
\equiv R_{m1}+R_{m2}+R_{m3}+R_{m4}, 
\ \ \mbox { for } m\geq 3
\end{cases}
\end{equation}
with $a_i=\frac{i+3/2}{m+1}$. The proof of the Lemma is complete.
\cbdu

We are now ready to prove Theorem \ref{Mthm:hm} by induction.

\medskip

{\bf{Proof of Theorem \ref{Mthm:hm}:}}
When $m=0$, inequality (\ref{dec:hm-u}) has been established in Lemma \ref{le:udecay} Section \ref{sec:3} and inequality (\ref{dec:hm1-d}) is proved in Theorem \ref{Mthm} (see \cite{DQS}).  Assume by induction that, for $k=0, 1, 2, . .., m-1$,
\begin{equation}\label{ass:u}
\|D^{k}u\|_{L^2(\R^3)}^2\leq C_K(1+t)^{-(k+3/2)},
\end{equation}
\begin{equation}\label{ass:d}
\|D^{k}d\|_{L^2(\R^3)}^2\leq C_k(1+t)^{-(k+3/2)}.
\end{equation}
We then apply the Fourier splitting method to the energy estimate (\ref{ineq:hm-estimate}). Lemma \ref{le:u-infty-decay}, Lemma \ref{le:d-infty} and inequality (\ref{dec:d2-infty}) in Lemma \ref{le:infty-interp} imply
\bg\notag
\|u\|_{L^\infty}^2+\|d-w_0\|_{L^\infty}^2+\|D^2 d\|_{L^\infty}^2\leq C(1+t)^{-1}.
\ed
Lemma \ref{le:d-infty} also implies 
\bg\notag
\|d-w_0\|_{L^\infty}\leq C(1+t)^{-1}.
\ed
The induction hypothesis (\ref{ass:d}) implies 
\bg\notag
\|D^{m-1}d\|_{L^2}^2\leq C_m(1+t)^{-(m+1/2)}.
\ed
By (\ref{ineq:hm-estimate}) we have 
\begin{align}\label{four:hm-estimate}
&\frac{d}{dt}\int_{\R^3}|D^mu|^2+|D^{m}d|^2dx+\int_{\R^3}|D^{m+1}u|^2+|D^{m+1}d|^2dx\\
&\leq C(1+t)^{-1}\int_{\R^3}|D^mu|^2+|D^{m}d|^2dx\notag\\
&+C(1+t)^{-(\frac{1}{2}+m+2)}+CR_m\notag.
\end{align}
For $m=1,2$, $R_m=0$.
For $m\geq 3$, we proceed to estimate each term $R_{mi}$ on the right hand side of (\ref{R}) as follows. The induction hypothesis (\ref{ass:u}) and the decay estimate (\ref{eq:udecay}) in Lemma \ref{le:udecay} yield
\begin{align}\label{R1}
R_{m1}&=\sum_{i=1}^{m/2}\|u\|_{L^2}^{2}\|D^{m-i}u\|_{L^2}^{2/(1-a_i)}\leq C_m(1+t)^{-\frac{3}{2}-\frac{m-i+3/2}{1-a_i}}\\
&\leq C_m(1+t)^{-\frac{3}{2}-(m+1)\frac{m-i+3/2}{m-i-1/2}}
\leq C_m(1+t)^{-(\frac{3}{2}+m+1)}\notag.
\end{align}
The hypothesis (\ref{ass:u}) and the decay estimate (\ref{ddecay}) in Theorem \ref{Mthm} yield
\begin{align}\label{R3}
R_{m3}&=\sum_{i=1}^{m-1}\|d-w_0\|_{L^2}^2\|D^{m-i}u\|_{L^2}^{\frac{2}{1-a_{i}}}
\leq C_m(1+t)^{-\frac{3}{2}-\frac{1/2+m-i+1}{1-a_{i}}}\\
&\leq C_m(1+t)^{-\frac{3}{2}-(m+1)\frac{m-i+3/2}{m-i-1/2}}
\leq C_m(1+t)^{-(\frac{3}{2}+m+1)}\notag.
\end{align}
The hypothesis (\ref{ass:d}) and the decay estimate (\ref{ddecay}) in Theorem \ref{Mthm} yield
\begin{align}\label{R2}
R_{m2}&=\sum_{i=2}^{(m-2)/2}\|d-w_0\|_{L^2}^{2}\|D^{m-i+1}d\|_{L^2}^{\frac{2}{1-a_{i+1}}}\\
&\leq C_m(1+t)^{-\frac{3}{2}-\frac{m-i+3/2+1}{1-a_{i+1}}}\notag\\
&\leq C_m(1+t)^{-\frac{3}{2}-(m+1)\frac{m-i+5/2}{m-i-1/2}}
\leq C_m(1+t)^{-(\frac{3}{2}+m+1)}\notag.
\end{align}
Similarly as to obtain (\ref{R2}), we have 
\begin{equation}\label{R4}
R_{m4}=\sum_{i=1}^{(m-1)/2}\|d-w_0\|_{L^2}^{2}\|D^{m-i+1}d\|_{L^2}^{\frac{2}{1-a_{i}}}\leq C_m(1+t)^{-(\frac{3}{2}+m+1)}.
\end{equation}
Combing the inequalities (\ref{four:hm-estimate})-(\ref{R4}) yields 
\begin{align}\label{four:energy}
&\frac{d}{dt}\int_{\R^3}|D^mu|^2+|D^{m}d|^2dx+\int_{\R^3}|D^{m+1}u|^2+|D^{m+1}d|^2dx\\
&\leq C_m(1+t)^{-1}\int_{\R^3}|D^mu|^2+|D^{m}d|^2dx+C_m(1+t)^{-(m+5/2)}\notag.
\end{align}
Applying the Fourier splitting method to the inequality (\ref{four:energy}) yields 
\begin{equation}\label{dec:udm}
\int_{\R^3}|D^mu|^2+|D^{m}d|^2dx\leq C_m(1+t)^{-(m+3/2)},
\end{equation}
for all $m\geq 1$. This proves inequalities (\ref{dec:hm-u}) and (\ref{dec:hm1-d}) for all $m\geq 1$.

Apply the following Gagliardo-Nirenberg interpolation inequality
\begin{equation}\notag
\|D^mu\|_{L^\infty(\R^3)}\leq C\|D^{m+2}u\|_{L^2(\R^3)}^{a}\|u\|_{L^2(\R^3)}^{1-a}
\end{equation}
with $a=\frac{2m+3}{2(m+2)}$. Combining (\ref{eq:udecay}) (in Lemma \ref{le:udecay}), (\ref{dec:udm}) and the last inequality yields
\begin{equation}\notag
\|D^mu\|_{L^\infty(\R^3)}\leq C_m(1+t)^{-(3/2+m+1)\frac{a}{2}-\frac{3a}{4}}=C_m(1+t)^{-\frac{m+3}{2}}. 
\end{equation}
This proves inequality (\ref{dec:infty-u}). Inequality (\ref{dec:infty-d}) can be proved similarly.
This completes the proof of Theorem \ref{Mthm:hm}.


{}

\end{document}